\documentstyle{czjphys209} % To be LaTeX'ed with czjphys209.sty

\begin{document}
\title{INHOMOGENEOUS MULTIPARAMETER JORDANIAN QUANTUM GROUPS 
       BY CONTRACTION}
\authori{Roger J. McDermott and Deepak Parashar}
\addressi{School of Computer and Mathematical Sciences,
The Robert Gordon University, \\ St. Andrew Street, Aberdeen AB25 1HG, U.K.\\
(e-mail: rm@scms.rgu.ac.uk, deeps@scms.rgu.ac.uk)}
\authorii{}
\addressii{}
\authoriii{}  % If only two authors
\addressiii{} % If only two authors
%\addressiii{Third address\\[1em]    %   
%{\normalsize\sc Fourth author} \\   %  If more than 3 authors
%\medskip Fourth address}            % 
\headtitle{Inhomogeneous Multiparameter Jordanian Quantum Groups \ldots}
\headauthor{Roger J. McDermott and Deepak Parashar} %use "et al." for more than 3 authors
\specialhead{R. J. McDermott: Jordanian Quantum Groups \ldots}
%%%%%%%%%%%%%%  FOR EDITORIAL USE ONLY!!!  %%%%%%%%%%%%%%%
\evidence{A}
\daterec{XXX}    %;\\ final version }
\cislo{0}  \year{1999}
\setcounter{page}{1}
\pagesfromto{000--000}
%\makefirsttitle
%%%%%%%%%%%%%%%%%%%%%%%%%%%%%%%%%%%%%%%%%%%%%%%%%%%%%%%%%%
\maketitle

\begin{abstract}
It is known that the inhomogeneous quantum group $IGL_{q,r}(2)$ can be constructed 
as a quotient of the multiparameter $q$--deformation of $GL(3)$. We show that a 
similar result holds for the inhomogeneous Jordanian deformation and exhibit its 
Hopf structure. 

\end{abstract}

\section{Introduction}

It is well--known \cite{[1]} that, analogous to the classical group--theoretical method, the $q$--deformation of $IGL(2)$ can be constructed by factoring out a certain two--sided Hopf ideal from the multiparameter $q$--deformation of $GL(3)$. This is an interesting procedure, allowing, for example, the construction of a differential calculus on the quantum plane by a reduction of the differential calculus on the quantum group. In this paper, we apply the same construction to the Jordanian deformation. The multiparameter Jordanian deformation of $GL(3)$ is first produced by a contraction from the corresponding $q$--deformation and this is then used to construct the inhomogeneous group by factorisation. The Hopf -structure of $IGL_J(2)$ is given explicitly and we show that it is possible to derive from this a coaction of a modified version of $GL_J(2)$ on the Jordanian quantum plane.\\

Note: In this paper, we denote  $q$--deformed structures using the (multiparameter) subscript $Q$ and structures that have been contracted to the Jordanian form are written with a subscript $J$ (e.g $GL_Q(3)$ and $GL_J(3)$). \\

\section{The R--matrix for GL$_Q$(2)}
Following Aschieri and Castellani [1], the $R$--matrix for $GL_Q(3)$ (where  $Q = \{ r,s,p,q \}$) can be written as 
\begin{equation}
R_Q(3)= 
\left( 
\begin{array}{cccc} 
r & & & \\ & S^{-1}& & \\ & \Lambda & S & \\ & & & R_Q(2)
\end{array}
\right)
\end{equation}
where
$S = \left( \begin{array}{cc} p & 0 \\ 0 & q \end{array}   \right)$, 
$\Lambda = \left(  \begin{array}{cc} r-r^{-1} & 0 \\ 0 & r-r^{-1} \end{array}   \right)$
and
\begin{equation}
R_Q(2)= \left( 
\begin{array}{cccc} r & & & \\ & s & & \\ & r-r^{-1} & s^{-1} & \\ & & & r \end{array}
\right)
\end{equation}
The matrix indices of $R_Q(3)$ run, in order, through the set 
$(11)$, $(12)$, $(13)$, $(21)$, $(31)$, $(22)$, $(23)$, $(32)$,
$(33)$.
This numbering system is chosen to clearly show the embedding of the
$R_Q(2)$ matrix in the $R_Q(3)$ matrix which, in turn, allows the Hopf
structure of larger algebra to be analysed in terms of the simpler one.
The Hopf structure of $GL_Q(3)$ is given by the $RTT$ relations with
$T$--matrix
\begin{equation}
{\cal T} = 
\left( 
\begin{array}{ccc} f & \theta & \phi \\ x & a & b \\ y & c & d \end{array} 
\right)
\end{equation}
and the multiparameter inhomogeneous $q$--deformation $IGL_Q(2)$ is the quantum homogeneous space 
\begin{equation}
IGL_Q(2) = GL_Q(3)/H
\end{equation}
where $H$ is the two--sided Hopf ideal generated by the $T$--matrix elements $\{ \theta , \phi \}$.

\section{The Contraction Procedure}
The $R$--matrix of the Jordanian (or $h$--deformation) can be viewed as a singular limit of a similarity transformation on the $q$--deformation $R$--matrix \cite{[2]}\cite{[3]}. Let $g(\eta)$ be a matrix dependent on a contraction parameter $\eta $ which is itself a function of one of the deformation parameters of the $q$--deformed algebra. This can be used to define a transformed $q$--deformed $R$--matrix
\begin{equation}
\tilde{R}_J = (g^{-1} \otimes g^{-1})R_Q(g \otimes g)
\end{equation}
The $R$--matrix of the Jordanian deformation is then obtained by taking a limiting value of the parameter $\eta$. Even though the contraction parameter $\eta$ is undefined in this limit, the new $R$--matrix is finite and gives rise to a new quantum group structure through the $RTT$--relations. For example, in the contraction process which takes $GL_q(2)$ to $GL_h(2)$, the contraction matrix is 
\begin{equation}
g(\eta) = \left( \begin{array}{cc} 1 & 0 \\ \eta & 1 \end{array} \right) \label{contracmat}
\end{equation}
where  $\eta = \frac{h}{1-q}$ with $h$ a new free parameter.\\

It has been shown by Alishahiha \cite{[3]} that, in the extension of this procedure to the construction of $GL_J(3)$, there are essentially two choices of contraction matrix. The first has been used in a number of papers, e.g. by Quesne \cite{[4]} and takes the form
\begin{equation}
G^{\prime} = \left( \begin{array}{ccc} 1 & 0 & 0 \\ 0 & 1 & 0 \\ \eta & 0 & 1 \end{array}
\right)
\end{equation}

There is, however, a second choice (also mentioned in \cite{[3]} but not pursued there since it gives trivial results for the single--parameter $q$--deformation) 
\begin{equation}
G = \left( \begin{array}{cc} 1 & 0 \\ 0 & g  \end{array} \right)
\end{equation}
where $g$ is the $2 \times 2$ contraction matrix
\begin{equation}
g(\eta) = \left( \begin{array}{cc} 1 & 0 \\ \eta & 1 \end{array} \right) 
\end{equation}
with $\eta = \frac{r}{1-q}$. In this present work, we take $G$ as our contraction matrix because, unlike the matrix $G^{\prime}$, after contraction it allows a non--trivial embedding of $R_J(2)$ in $R_J(3)$ in a manner similar to the $q$--deformed case. It is then possible to perform the quotient construction for the inhomogeneous quantum group.\\

If the similarity transformation is made using the matrix $G$, we obtain
\begin{eqnarray}
R_J(3)  &=&  \lim_{r \rightarrow 1} 
\left( \begin{array}{cccc} 
r & & & \\
  & g^{-1}S^{-1}g & & \\
  & \Lambda & g^{-1}Sg & \\
  & & & (g^{-1} \otimes g^{-1})R_Q(g \otimes g)
\end{array} \right) \\
 & = & 
\left( \begin{array}{cccc} 
1 & & & \\
  & K^{-1} & & \\
  &  & K & \\
  & & & R_J(2)
\end{array} \right)
\end{eqnarray} 
where 
$K$ is the matrix 
$\left( \begin{array}{cc} p & 0 \\ k & p \end{array}
\right)$
and $R_J(2)$ is the $R$--matrix for the multiparameter Jordanian deformation of $GL(2)$
\begin{equation}
\left(
\begin{array}{rrrr} 1 & & & \\ m & 1 & & \\ -m & 0 & 1 & \\ mn & n & -n & 1 \end{array}
\right)
\end{equation}
The free parameters $\{m,n,k\}$ appear as limits in the contraction process while the parameter $\{ p \}$ survives the contraction process. The result is a four parameter Jordanian defomation of $GL(3)$.

\section{Multiparameter Jordanian Deformation of GL(3)}

We denote the $T$--matrix for the Jordanian deformation by 
\begin{equation}
\cal{T} = \left(
		\begin{array}{ccc} f & \theta & \phi \\ x & a & b \\ y & c & d \end{array} 
	  \right)
	= \left(
	        \begin{array}{cc} f & \Theta  \\ X & T \end{array}
	  \right)
\end{equation}
where $T = \left( \begin{array}{cc} a & b  \\ c & d \end{array} \right)$, 
$X = \left( \begin{array}{c} x  \\ y \end{array} \right)$, and  $\Theta = \left( \begin{array}{cc} \theta , \phi \end{array} \right)$. 
\subsection{Algebra Relations}
The algebra structure of the quantum group is obtained through the $RTT$--procedure. For the commutation relations between the elements of the matrix $T$, we have the usual relations between the generators of the multiparametric Jordanian deformation of $GL(2)$: 
\begin{equation}
\begin{array}{lll}
[a,b] =  nb^2  & [a,c] = m(\delta - a^2)   & [a,d] = nbd-mba  \\  
{[}{b},{d}{]} = -mb^2 & [b,c] = -mba -ndb  & [c,d] = n(d^2 - \delta ) 
\end{array}
\end{equation}
where $\delta$ is the quantum determinant of the submatrix $T$
\begin{equation}
\delta = ad - bc -nbd
\end{equation}
with commutation relation 
\begin{equation}
\begin{array}{llll}
[\delta, a] = (m-n) \delta b & [\delta, b] = 0  & [\delta , c] = (m-n) (\delta d - a \delta) &
[\delta, d] = (n-m) \delta b 
\end{array}
\end{equation}
Thus $\delta$ is not central in the (sub--Hopf) algebra generated by elements of $T$ unless $m=n$.\\

The relations between $T$ and $f$ are given by
\begin{equation}
\begin{array}{cccc}
[a,f] = \frac{k}{p}fb & [b,f] = 0 & [c,f] = \frac{k}{p}(fd -af) & [d,f] = -\frac{k}{p}bf
\end{array}
\end{equation}
those between $T$ and $X$ are 
\begin{equation}
\begin{array}{ll}
[a,x]_p = kxb & [b,x]_p = 0 \\ {[}c,x]_p = kxd + max & [d,x]_p = mbx \\
{[}{a},{y}{]}_p = kyb - max & [b,y]_p = -mbx \\ {[}c,y]_p = kyd + ncx -nay - mnax & 
[d,y]_p = ndx - nby -mnbx \\
\delta x = p^2 x \delta &  \delta y =p^2 y \delta +(n-m) \delta x 
\end{array}
\end{equation}
while those between  $f$ and $X$ give
\begin{equation}
\begin{array}{cc}
[f,x]_p =0 & [f,y]_p = -kxf
\end{array}
\end{equation}
The commutation relations between the elements of $X$ are the usual relations for the Jordanian quantum plane $C_J(2)$:
\begin{equation}
[x,y] = -mx^2
\end{equation}
There are also similar commutation relations between the elements of $T$, $f$ and $\Theta$, as well as cross-relations between $X$ and $\Theta$.

\subsection{Coalgebra Relations and Antipode}

The coalgebraic structure of the Hopf algebra is the usual one:
\begin{equation}
\Delta({\cal T}) = {\cal T} \dot{\otimes} {\cal T} \hspace{10mm} \epsilon({\cal T}) = I_3
\end{equation}
with antipode
\begin{equation}
S({\cal T}) = 
\left(
\begin{array}{cc}    e      & -e\Theta T^{-1} \\
		  -T^{-1}Xe & T^{-1}Xe \Theta T^{-1} + T^{-1}
\end{array}
\right)
\end{equation}
where we append to the algebra, the element $e = (f- \Theta T^{-1}X)^{-1}$. In terms of these elements, the quantum determinant of the $T$-matrix ${\cal T}$ is 
\begin{equation}
{\cal D} = \det ( {\cal T} ) = {e}^{-1} \delta
\end{equation}
and so, in the usual way, we can add $\xi = {\cal D}^{-1}$ to the algebra to obtain the full Hopf algebra.

\section{The Inhomogeneous Multiparameter Jordanian Quantum Group IGL$_J$(2)}

We define $H$ to be the space of all monomials containing at least one element of $\Theta$. It is straightforward to prove the following :

\begin{enumerate}
\item $H$ is a two--sided ideal in $GL_J(3)$.

\item $H$ is a co--ideal i.e. $\Delta (H) \subseteq H\otimes GL_J(3) + GL_J(3) \otimes H$ and $\epsilon (H) = 0$.

\item $S(H) \subseteq H$.
\end{enumerate}

Thus $H$ is a two--sided Hopf ideal and so we can define a canonical projection from $GL_J(3)$ to the quotient space $GL_J(3)/H$ which respects the Hopf--algebraic structure (i.e. the $RTT$--relations). Consequently the quotient is a Hopf algebra which we denote $IGL_J(2)$.\\

The algebra sector for this quantum group has commutation relations formally obtained from $GL_J(3)$ by setting the generator set $\Theta = 0$ and this gives rise to the commutation relations explicitly detailed in the previous section. The $T$--matrix for the coalgebra is given by
\begin{equation}
{\cal T} =  \left( \begin{array}{cc} f & 0  \\ X & T \end{array} \right)
\end{equation}
which gives the coproduct
\begin{equation}
\Delta({\cal T}) = {\cal T} {\dot{\otimes}} {\cal T} =
\left( \begin{array}{cc} 
f\otimes f & 0 \\ 
T\dot{\otimes} X + X \dot{\otimes} f & T \dot{\otimes} T
\end{array} \right)
\end{equation}
counit $\epsilon ({\cal T}) = I_{3}$ and antipode
\begin{equation}
S({\cal T}) = 
\left( \begin{array}{cc} f^{-1} & 0 \\ -T^{-1}Xf^{-1} & T^{-1} \end{array} \right)
\end{equation}
The quantum determinant ${\cal D} = f \delta$ is group--like but, since $f$ is not central, it cannot be made simultaneously central with $\delta$ unless the whole algebraic structure collapses to a trivial extension of the single--parameter Jordanian deformation of $GL(2)$. This is analagous to the situation in the $q$--deformed case.\\

This procedure also shows that it is possible to view the Jordanian
quantum plane, $C_J(2)$,  as the quantum homogeneous space
$IGL_J(2)/GL_J(2)^*$ where $GL_J(2)^*$ is the Hopf algebra formed by
appropriately appending the ``dilatation element'' $f$ to $GL_J(2)$. The
comultiplication in $IGL_J(2)$ can then be viewed as a coaction of the
quantum group $GL_J(2)^*$ on the quantum plane $C_J(2)$ generated by the
elements $X$. However, unlike the usual case, there is a non--trivial
braiding between the elements of the quantum group and quantum plane.

\section{Conclusion}

We have shown that it is possible to construct the inhomogeneous Jordanian deformation $IGL_J(2)$ as a quotient group by factoring out a Hopf ideal from $GL_J(3)$. It would be of interest to construct the differential calculus on the Jordanian quantum plane by a reduction of the bicovariant differential caculus on $GL_J(3)$. This would be allow the investigation of physical models with $GL_J(N)$ symmetry similar to that of Cho {\it et al} \cite{[5]} and Madore and Steinacker \cite{[6]}. Work on this problem is underway.

\poslednisuda % uncomment if your paper has even number of pages
\end{document}